\documentclass[11pt,a4paper]{article}

\usepackage[T1]{fontenc}
\usepackage[sc]{mathpazo}
\usepackage{microtype}
\usepackage{amsmath,amssymb,amsthm,mathtools}
\usepackage[margin=2.5cm]{geometry}
\usepackage{booktabs,array}
\usepackage{graphicx}
\usepackage{xcolor}
\usepackage{enumitem}
\usepackage[colorlinks=true,linkcolor=blue!50!black,citecolor=blue!50!black,urlcolor=blue!50!black]{hyperref}

\newtheorem{theorem}{Theorem}[section]
\newtheorem{lemma}[theorem]{Lemma}
\newtheorem{proposition}[theorem]{Proposition}
\newtheorem{corollary}[theorem]{Corollary}
\theoremstyle{definition}
\newtheorem{definition}[theorem]{Definition}
\newtheorem{example}[theorem]{Example}
\theoremstyle{remark}
\newtheorem{remark}[theorem]{Remark}

\DeclarePairedDelimiter{\abs}{\lvert}{\rvert}
\DeclarePairedDelimiter{\norm}{\lVert}{\rVert}
\newcommand{\R}{\mathbb R}
\newcommand{\Bu}{\mathbb B}
\newcommand{\A}{\mathcal A}
\newcommand{\Hq}{\mathbb H}
\newcommand{\dbar}{\overline{\partial}}
\newcommand{\Cst}{C_*}
\newcommand{\Hs}{\mathsf H}
\newcommand{\diag}{\operatorname{diag}}
\providecommand{\ip}[2]{\langle #1,#2\rangle}
\providecommand{\Reop}{\operatorname{Re}}
\title{\textbf{Estimates of the Schwarz constant for monogenic functions in $\mathbb R^3$}}
\author{Hung Manh Nguyen$^{(1)}$}
\date{\small $^{(1)}$ University of Transport and Communications, Hanoi, Vietnam\\ Email: hung.manh.nguyen@utc.edu.vn}

\begin{document}
\maketitle

\begin{abstract}
The Schwarz lemma bounds the differential at the origin of a holomorphic self-map of the unit disc fixing $0$: $\abs{f'(0)}\le1$. For monogenic maps in $\mathbb R^3$, the optimal constant is not known. Let $\mathcal M$ be the class of monogenic functions $f$ from the unit ball of $\mathbb R^3$ to the reduced quaternions with $\abs f\le1$ and $f(0)=0$, and let $\Cst=\sup_{\mathcal M}\abs{f'(0)}$, where $f'=\partial_0f$ is the hypercomplex derivative. We show that $\frac{80}{61}\le\Cst\le1.35784$. Thus the classical bound $1$ fails, while $\Cst$ stays below the sharp constant $\frac32$ for harmonic maps of the ball into itself. The estimates rest on a reduction, by rotation covariance and averaging about the $x_0$-axis, to a convex extremal problem in one variable for zonal harmonic functions. Lower bounds come from restricting this problem to polynomials, and upper bounds from its dual formulation in terms of vector fields on the sphere.
\end{abstract}


\section{Introduction}\label{sec:intro}

Let $f$ be holomorphic in the unit disc with $\abs{f}\le1$ and $f(0)=0$. Schwarz's lemma gives $\abs{f'(0)}\le1$ and $\abs{f(z)}\le\abs{z}$. The proof uses three tools: the factorization $f=zg$, the maximum principle for $g$, and the Cauchy estimate with a kernel of constant modulus on the circle, which gives $\abs{f'(0)}\le1$ directly. Geometrically, the differential of a holomorphic map is a similarity, and the lemma bounds
its dilation factor $\abs{f'(0)}$.

In $\R^3$, monogenic functions are a natural analogue of holomorphic functions: they are the null solutions of the Cauchy--Riemann operator, with values in the reduced quaternions \cite{BDS82,GS97,SW68}. Their derivative is the hypercomplex derivative $f'=\frac12\partial f = \partial_0f$ \cite{GM99}. Their differential is not a similarity: at a point where it is non-degenerate, it maps infinitesimal spheres onto ellipsoids whose longest semi-axis equals the sum of the other two \cite{GNM12}. Consequently, there is no factorization $f=xg$, and the classical proof breaks down. Yang and Qian \cite{YQ06} prove a sharp Schwarz lemma for monogenic functions outside the unit ball, with the Cauchy kernel as extremal, and remark that the interior version is false.

Inside the ball, Schwarz-type inequalities with explicit but non-optimal constants were obtained by Zhang \cite{Zh14} for Clifford-valued monogenic functions vanishing at a point $a$, and extended by Xu, Yu and Huo \cite{XYH25} to harmonic Clifford-valued functions. At $a=0$, they give $\abs{f(x)}\le\abs x/(\sqrt[3]2-1)\approx3.85\,\abs x$ in $\R^3$, which is
rather large compared with the constant $1$ of the classical lemma. To our knowledge, the sharp theory that does exist is for harmonic functions, where the monogenic structure plays no role. For real harmonic $u$ on the unit ball of $\R^n$ with $\abs u\le1$, the best bound at the center is $\abs{\nabla u(0)}\le2\abs{\mathbb B_{n-1}}/\abs{\mathbb B_n}$ \cite{ABR01}, where $\abs{\mathbb B_n}$ is the volume of the unit ball of $\R^n$. This constant equals $\frac32$ for $n=3$. Kalaj \cite{Ka} gives the sharp constant for the differential at the origin of harmonic maps into $\R^m$ in Hardy norms. Liu \cite{Li22} proves the sharp Schwarz--Pick inequality for real-valued harmonic functions at all points, with the curious feature that for $n=3$ the constant away from the center is $8/(3\sqrt3)$ rather than $\frac32$. Forstneri\v{c} and Kalaj \cite{FK} obtain a Schwarz--Pick lemma for harmonic maps of the disc into a ball which are conformal at a point, with equality exactly for conformal diffeomorphisms. By contrast, the differential of a monogenic map in $\R^3$ is never conformal unless it vanishes.

The Schwarz lemma is a starting point of broader theories. Its invariant form, the Schwarz--Pick lemma, says that holomorphic self-maps of the disc contract the hyperbolic metric, and it is closely related to Schur analysis, the study of functions analytic and contractive in the disc (Schur functions) \cite{A2001}. D.~Alpay and his collaborators carried these ideas over to slice regular functions of a quaternionic variable: the slice regular self-maps of the ball are identified as contractive multipliers of the Hardy space, and a Nevanlinna--Pick interpolation theory is developed, from which the Schwarz--Pick lemma follows \cite{ABCS15}. Slice regular functions carry a product structure that monogenic functions of three real variables lack. Whether the methods of Schur analysis can be brought to bear on monogenic functions, starting with a better understanding of the Schwarz lemma in this setting, remains open.

We ask for the best constant in the Schwarz inequality for the hypercomplex derivative. Let
\[
\mathcal M=\{f:\Bu\to\A\ \text{monogenic},\ \abs f\le1,\ f(0)=0\},\qquad \Cst=\sup_{f\in\mathcal M}\abs{f'(0)},
\]
where $\Bu$ is the open unit ball of $\R^3$ and $\A$ the reduced quaternions. The classical statement would be $\Cst=1$. The contributions of this paper are as follows.
\begin{enumerate}[label=(\roman*),itemsep=1pt]
\item An explicit cubic shows $\Cst\ge\frac{11}{10}$, and a Cauchy-kernel estimate gives $\Cst\le C_K=1+\ln(2+\sqrt3)/(2\sqrt3)$, approximately $1.38017$. The sharp constant for all harmonic maps is $\frac32$, so monogenicity strictly lowers it (Section~\ref{sec:cauchy}).
\item Rotation covariance shows that $\Cst$ is simultaneously the supremum of $\abs{f'(0)}$, of the operator norm of the Jacobian, and of a single entry $\partial_0f_0(0)$. Averaging about $x_0$-axis reduces the extremal problem to a convex problem for zonal harmonic functions (Section~\ref{sec:red}), which has a primal and a dual formulation (Section~\ref{sec:duality}).
\item We prove the bounds $\frac{80}{61}\le\Cst\le 1.35784$ (Example~\ref{ex:lower} and Corollary~\ref{cor:U}); the upper bound comes from an explicit dual field, while the lower bound comes from a polynomial of degree $6$.
\end{enumerate}

Section~\ref{sec:prelim} presents notations and recalls the gradient representation of monogenic functions. Section~\ref{sec:bounds} gives two elementary bounds for $\Cst$. Section~\ref{sec:extremal} reduces the extremal problem by symmetry to a one-dimensional convex problem, studies its primal and dual formulations, and derives sharper lower and upper bounds and the comparison with harmonic maps. Section~\ref{sec:concl} summarizes the results and discusses open problems.

\section{Preliminaries}\label{sec:prelim}

\subsection{Notation and conventions}\label{sec:notation}

\paragraph{Quaternions.}
Let $\Hq$ be the {real associative algebra with basis $e_0,e_1,e_2,e_3$, where $e_0=1$ is the unit and
\[
e_1e_2=e_3;\quad e_ie_j+e_je_i=-2\delta_{ij},\quad i,j=1,2,3.
\]
For every $p=\sum_{i=0}^{3}p_ie_i\in\Hq$ and $q=\sum_{i=0}^{3}q_ie_i\in\Hq$ with real $p_i,q_i$, we define the conjugate of $p$, the norm of $p$, and the inner product of $p$ and $q$ as follows:
\[
\overline{p}=p_0e_0-p_1e_1-p_2e_2-p_3e_3,\qquad \abs p=\ip pp^{1/2},\qquad \ip pq=\sum_{i=0}^3p_iq_i.
\]
We use the following standard facts: 
$$\overline{pq}=\overline{q}\,\overline{p};\qquad p\,\overline{p}=\overline{p}\,p=\abs p^2;\qquad \ip pq=\Reop(p\overline{q}),$$
where $\Reop(p)=p_0$. Consequently the norm is multiplicative,
\[
\abs{pq}^2=pq\,\overline{pq}=p\,(q\,\overline{q})\,\overline{p}=\abs q^2\,p\,\overline{p}=\abs p^2\abs q^2 ,
\]
and conjugation is an isometry: $\abs{\overline{p}}=\abs p$ and $\ip{\overline{p}}{\overline{q}}=\ip pq$.}

\paragraph{Reduced quaternions and $\R^3$.}
The space of reduced quaternions is $\A=\operatorname{span}_\R\{e_0,e_1,e_2\}\subset\Hq$. It is a real vector subspace, but not a subalgebra ($e_1e_2=e_3\notin\A$). We identify $\R^3$ with $\A$ by
\[
x=(x_0,x_1,x_2)\ \longleftrightarrow\ x=\sum_{i=0}^{2}x_ie_i=x_0e_0+x_1e_1+x_2e_2 ,
\]
so that $\abs x$ and $\ip xy$ are the Euclidean norm and inner product of $\R^3$, and the vector $e_0=(1,0,0)$ is the first unit vector. Likewise, a column vector in $\R^3$ is identified with the element of $\A$ having the same coordinates. For $x\in\R^3$ we put
\[
r=\abs x,\qquad \rho=\sqrt{x_1^2+x_2^2},\qquad r^2=x_0^2+\rho^2,\qquad x_0=r\cos\vartheta\ \ (x\ne0),
\]
where $\vartheta$ is the polar angle between $x$ and the $x_0$-axis.

\paragraph{Derivatives and Cauchy--Riemann operators.}
{We denote $\partial_i=\partial/\partial x_i$ for $i=0,1,2$ and for $w\in\R^3$, $\partial_w=\sum_{i=0}^2w_i\partial_i$. Acting on $\Hq$-valued functions from the left, the Cauchy--Riemann} operator and its conjugate are
\[
\dbar=\sum_{i=0}^{2}e_i\partial_i=\partial_0+e_1\partial_1+e_2\partial_2,\qquad
\partial=\sum_{i=0}^{2}\overline{e_i}\,\partial_i=\partial_0-e_1\partial_1-e_2\partial_2 .
\]
Obviously, we have $\partial+\dbar=2\partial_0$, and the decomposition of the Laplace operator 
$$\dbar\,\partial=\partial\,\dbar=\Delta:=\partial_0^2+\partial_1^2+\partial_2^2.$$
A function $f\in C^1(\Omega;\Hq)$ on an open set $\Omega\subset\R^3$ is (left) monogenic if $\dbar f=0$ in $\Omega$.

For a real-valued $h\in C^1$, we gets
\[
\dbar h=\sum_{i=0}^{2}(\partial_ih)\,e_i\in\A,\qquad \partial h=\overline{\dbar h}=\partial_0h-\partial_1h\,e_1-\partial_2h\,e_2 .
\]
Thus $\dbar h$ is the gradient of $h$ written as an element of $\A$, and $\abs{\dbar h}=\abs{\partial h}=\Big(\sum_{i=0}^2(\partial_ih)^2\Big)^{1/2}$.
The Hessian matrix of a real-valued $h\in C^2$ is defined as
$$\Hs_h(x)=\Big(\partial_i\partial_jh(x)\Big)_{i,j=0}^{2}.$$
It is a real symmetric $3\times3$ matrix, and $\operatorname{tr}\Hs_h=\Delta h$. For a real $3\times3$ matrix $M$ and a column vector $w\in\R^3$, $Mw\in\R^3\cong\A$ denotes the usual matrix--vector product.

\paragraph{Differential and hypercomplex derivative.}
For $f=\sum_{i=0}^2f_ie_i\in C^1(\Omega;\A)$ and $x\in\Omega$, the differential of $f$ at $x$ is the linear map 
$$df_x:\R^3\to\A\, ,\qquad df_x(w)=\partial_wf(x)=\sum_{j=0}^2w_j\,\partial_jf(x),$$
where $\partial_wf(x)$ is the directional derivative of $f$ at $x$ in the direction $w$. In the bases $e_0,e_1,e_2$, it is represented by the Jacobian matrix
\[
J_f(x)=\Big(\partial_jf_i(x)\Big)_{i,j=0}^{2}\, ,\qquad J_f=J_f(0),
\]
so that $J_f(x)w=df_x(w)=\partial_wf(x)$. We write $J_f(x)w$ when we regard $J_f(x)$ as the differential (or as a matrix), and $\partial_wf(x)$ when we think of a derivative in a given direction. The Hessian matrix is the Jacobian matrix of the gradient, $\Hs_h(x)=J_{\dbar h}(x)$, so that $\Hs_h(x)w=\dbar(\partial_wh)(x)$. For monogenic $f$, the hypercomplex derivative is $f'=\frac12\partial f = \partial_0f$ (see \cite{GM99}), i.e. the first column of $J_f$, thus $f'(0)=\partial_0f(0)$.

\paragraph{Balls, norms.}
For $z\in\R^3$ and $R>0$, $B(z,R)=\{x\in\R^3:\abs{x-z}<R\}$ is the open ball and $\overline B(z,R)$ is its closure. We put $\Bu=B(0,1)$, $S^2=\partial\Bu$, and let $\sigma = 4\pi$ be the surface measure of $S^2$.
Three different norms are used, and we keep them distinct:
\begin{itemize}[itemsep=1pt]
\item $\abs\cdot$ is the Euclidean norm on $\R^3\cong\A$ and on $\Hq$;
\item for a real $3\times3$ matrix $M$, $\norm M_{\mathrm{op}}=\max_{\abs w=1}\abs{Mw}$ is the operator norm; for a symmetric $M$, $\norm M_{\mathrm{op}}=\max_{\abs w=1}\abs{\ip{Mw}w}$ is the largest modulus of an eigenvalue;
\item for a bounded function $g$ on a set $E$, $\norm g_\infty=\sup_E\abs g$.
\end{itemize}

\subsection{Monogenic functions are conjugated gradients}

\begin{lemma}[Conjugate-gradient representation]\label{lem:grad}
Let $\Omega\subset\R^3$ be an open set and star-shaped with respect to $0$ and $f\in C^1(\Omega;\A)$.
\begin{enumerate}[label=(\alph*),itemsep=1pt]
\item $f$ is monogenic if and only if $f=\partial h$ for a real harmonic function $h$ on $\Omega$; $h$ is unique up to an additive constant.
\item If $f=\partial h$ is monogenic, then all components of $f$ are harmonic, $f'=\tfrac12\partial f=\partial_0f$, and for all $x\in\Omega$, $w\in\R^3$
\[
\begin{gathered}\abs{f(x)}=\abs{\partial h(x)},\qquad \partial_wf(x)=\partial(\partial_wh)(x)=\overline{\Hs_h(x)w},\\ 
J_f(x)=C\,\Hs_h(x),\quad C=\operatorname{diag}(1,-1,-1),\end{gathered}
\]
where $\Hs_h(x)$ is symmetric and trace-free. In particular 
$$
\begin{gathered}
f'(x)=\partial(\partial_0h)(x)=\overline{\Hs_h(x)e_0},\qquad \abs{\partial_wf(x)}=\abs{\Hs_h(x)w},\\
\qquad \norm{J_f(x)}_{\mathrm{op}}=\norm{\Hs_h(x)}_{\mathrm{op}},\qquad \partial_0f_0(x)=\partial_0^2h(x)=\ip{\Hs_h(x)e_0}{e_0},
\end{gathered}
$$ 
where $f_0=\Reop f$.
\end{enumerate}
\end{lemma}
\begin{proof}
(a) This is the classical fact that the conjugate of a monogenic function $f$ is the gradient of a harmonic function $h$, i.e. the components of $\overline{f}$ satisfy the generalised Cauchy--Riemann (M. Riesz) system (see Stein and Weiss \cite{SW68}; see also \cite{BDS82,GHS08,GS97}). We include a proof in the present notation, since all later computations are based on it.

Write $f=f_0+f_1e_1+f_2e_2$. We have
\begin{align*}
\dbar f&=\partial_0f+e_1\partial_1f+e_2\partial_2f\\
&=(\partial_0f_0-\partial_1f_1-\partial_2f_2)\,e_0+(\partial_0f_1+\partial_1f_0)\,e_1+(\partial_0f_2+\partial_2f_0)\,e_2+(\partial_1f_2-\partial_2f_1)\,e_3 .
\end{align*}
Put $F=\overline{f}=F_0+F_1e_1+F_2e_2$, i.e.\ $F_0=f_0$, $F_1=-f_1$, $F_2=-f_2$. Then $\dbar f=0$ is equivalent to the system
\begin{equation}\label{eq:riesz}
\sum_{i=0}^{2}\partial_iF_i=0,\qquad \partial_jF_i=\partial_iF_j\quad(i,j=0,1,2).
\end{equation}
If $\dbar f=0$, the second part of \eqref{eq:riesz} says that the $1$-form $\sum_{i=0}^2 F_i\,dx_i$ is closed. On the star-shaped domain $\Omega$, it is exact by the Poincar\'e lemma \cite{Sp65}, and
$$h(x)=\int_0^1\sum_{i=0}^2 F_i(tx)x_i\,dt$$ 
satisfies $\partial_ih=F_i$. Then $h\in C^2$ and $\Delta h=\sum_{i=0}^2\partial_iF_i=0$. In quaternionic form $\dbar h=F=\overline{f}$, i.e. $f=\partial h$. In addition, two such $h$ have the same gradient, so they differ by a constant if $\Omega$ is connected.
Conversely, if $h$ is real harmonic, hence $h\in C^\infty(\Omega)$, and $f=\partial h$. Then $\dbar f=\dbar\partial h=\Delta h=0$.

(b) The components of $f=\partial h$ are $\partial_0h,-\partial_1h,-\partial_2h$, which are harmonic since $h$ is. From $\dbar f=0$ we get $\frac12\partial f=\frac12(2\partial_0-\dbar)f=\partial_0f$. Next, $\abs f=\abs{\partial h}$ by definition. For $w\in\R^3$, the operators $\partial_w$ and $\partial$ have constant coefficients and commute, so $\partial_wf=\partial(\partial_wh)=\overline{\dbar(\partial_wh)}$, and 
\[
\dbar(\partial_wh)=\sum_{i=0}^2(\partial_i\partial_wh)\,e_i=\sum_{i=0}^2\sum_{j=0}^2 w_j\partial_i\partial_jh\,e_i=\Hs_hw.
\]
Hence $\partial_wf=\overline{\Hs_hw}$. Taking $w=e_j$ gives the $j$-th column of $J_f$, namely $\overline{\Hs_he_j}$. Since conjugation on $\A$ is $y\mapsto Cy$ in the bases $e_0,e_1,e_2$, this says $J_f=C\Hs_h$. The case $w=e_0$ gives $f'=\partial_0f=\overline{\Hs_he_0}$, and the $e_0$-component of $\overline{\Hs_he_0}$ is $(\Hs_he_0)_0=\partial_0^2h$, which gives $\partial_0f_0=\partial_0^2h=\ip{\Hs_he_0}{e_0}$. Conjugation is an isometry and $C$ is orthogonal, which gives the equalities of norms. $\Hs_h$ is symmetric because $h\in C^2$, and trace-free because $\Delta h=0$.
\end{proof}

All quantities below are therefore expressed through $\partial h$ and the Hessian $\Hs_h$. We also record the algebraic fact behind the ellipsoid property of M-conformal maps \cite{GNM12}.

\begin{corollary}[Singular values of the differential; cf.\ \cite{GNM12}]\label{cor:sing}
Let $f=\partial h$ be monogenic near $x$ and let $\lambda_1,\lambda_2,\lambda_3$ be the eigenvalues of $\Hs_h(x)$, then $\sum_{i=0}^2\lambda_i=0$. The singular values $s_1\ge s_2\ge s_3\ge0$ of $J_f(x)$ are the numbers $\abs{\lambda_i}$ in decreasing order, and they satisfy
\begin{equation}\label{eq:axis}
s_1=s_2+s_3 .
\end{equation}
In particular, if $J_f(x)\ne0$ then $s_3\le\frac12s_1<s_1$, meaning that $J_f(x)$ is not a similarity (it does not maps $S^2$ onto a sphere), and $\dfrac{s_1}{s_3}\ge2$ whenever $s_3>0$.
\end{corollary}
\begin{proof}
By Lemma~\ref{lem:grad}(b), $J_f(x)=C\,\Hs_h(x)$, thus $J_f(x)$ and $\Hs_h(x)$ have the same singular values. Since $\Hs_h(x)$ is symmetric, its singular values are $\abs{\lambda_i}$, and 
$$\sum_{i=0}^2\lambda_i=\operatorname{tr}\Hs_h=\Delta h=0.$$
Among three real numbers, two have the same weak sign (both $\ge0$ or both $\le0$); let us denote by $\lambda_a,\lambda_b$. The third is $\lambda_c=-(\lambda_a+\lambda_b)$, and $\abs{\lambda_c}=\abs{\lambda_a}+\abs{\lambda_b}\ge\abs{\lambda_a},\abs{\lambda_b}$. Hence $s_1=\abs{\lambda_c}=s_2+s_3$. Finally $s_2\ge s_3$ gives $s_1\ge2s_3$.
\end{proof}

\begin{lemma}[Maximum modulus; cf.\ \cite{BDS82,GHS08}]\label{lem:max}
Let $f:\overline\Bu\to\A$ be countinous on $\overline\Bu$ and monogenic in $\Bu$. Then 
$\abs{f(x)}\le\max_{y\in S^2}\abs {f(y)},\ x\in\overline\Bu.$ In particular, if $\abs f\le1$ on $S^2$, then $\abs f\le1$ in $\Bu$.
\end{lemma}

\subsection{Zonal harmonics and Appell polynomials}\label{sec:appell}

\paragraph{Legendre polynomials.}
The Legendre polynomials $P_n$ ($n\ge0$) are defined by the generating function
\begin{equation}\label{eq:genP}
(1-2ct+t^2)^{-1/2}=\sum_{n\ge0}P_n(c)\,t^n,\qquad \abs c\le1,\ \abs t<1 ,
\end{equation}
see e.g.\ \cite{Sz75}. In particular, some of the first Legendre polynomials are 
\[
\aligned
& P_0=1, \qquad\qquad P_1(c)=c, \qquad\qquad P_2(c)=\tfrac12(3c^2-1),\\
& P_3(c)=\dfrac12(5c^3-3c), \qquad P_4(c)=\dfrac18(35c^4-30c^2+3).
\endaligned
\]
We shall also use the following classical properties of the Legendre polynomials $P_n$:
\begin{itemize}
\item $P_n$ is a polynomial of degree $n$, has parity $P_n(-c) = (-1)^nP_n(c)$, and is normalized by $P_n(1)=1$;
\item the polynomials $P_n$ are orthogonal on $[-1, 1]$:
$$\int_{-1}^1P_n(c)P_m(c)\,dc=\frac2{2n+1}\delta_{nm};$$
\item $P_n$ satisfies the Legendre's differential equation
$$(1-c^2)P_n''(c) - 2cP_n'(c) + n(n+1)P_n(c) = 0.$$
\end{itemize}

\paragraph{Zonal harmonic polynomials.}
For $x\in\Bu$ and $\abs t<1$, we have
\[
1-2x_0t+r^2t^2=(1-x_0t)^2+\rho^2t^2>0 ,
\]
and we define $R_n(x)$ by the generating function
\begin{equation}\label{eq:genR}
\Phi(x,t):=(1-2x_0t+r^2t^2)^{-1/2}=\sum_{n\ge0}R_n(x)\,t^n.
\end{equation}
The series converges for $|t|<1$: for $x\neq0$ this follows from \eqref{eq:genP} with $c=x_0/r$ and $t$ replaced by $rt$ (since $|rt|<1$), and for $x=0$ we have $\Phi\equiv1$.

\begin{lemma}[cf.\ {\cite{ABR01}}]\label{lem:R}
For every $n\ge0$, the following properties hold:
\begin{enumerate}[label=(\roman*),itemsep=1pt]
\item $R_n$ is a homogeneous polynomial of degree $n$ in $x_0,x_1,x_2$, depending on $x_1,x_2$ only through $\rho^2=x_1^2+x_2^2$. Moreover,
$R_n(x)=r^nP_n\left(\dfrac{x_0}{r}\right)=r^nP_n(\cos\vartheta)$ for $x\ne0.$
\item $R_n(te_0)=t^n$ for $t\in(-1,1)$. In particular, $\partial_1R_n=\partial_2R_n=0$ on the $x_0$-axis.
\item $R_n$ is harmonic.
\item For every $n\ge0$, $\partial_0R_{n+1}=(n+1)R_n$.
\end{enumerate}
\end{lemma}
\begin{proof}
(i) Let $s=2x_0t-r^2t^2$. For $x\in\Bu$ and $|t|<\sqrt{2}-1$, then $|s|<1$ and we may expand 
$$(1-s)^{-1/2}=\sum_{k\ge0}\binom{2k}{k}4^{-k}s^k.$$ 
As $s^k=t^k(2x_0-r^2t)^k$, the coefficient of $t^n$ is a finite combination of terms $x_0^{k-j}(r^2)^{j}$ with $k+j=n$. Therefore $R_n$ is a homogeneous polynomial of degree $n$ in $x_0,\rho^2$ (note that $r^2=x_0^2+\rho^2$). For $x\ne0$, applying \eqref{eq:genP} with $c=x_0/r$ and $t$ replaced by $rt$ gives
$$\Phi(x,t)=\sum_{n\ge0} P_n\left(\dfrac{x_0}{r}\right)(rt)^n,$$ 
hence $R_n=r^nP_n\left(\dfrac{x_0}{r}\right)$.

(ii) On the $x_0$-axis, 
$$\Phi(\tau e_0,t)=(1-2\tau t+\tau^2t^2)^{-1/2}=(1-\tau t)^{-1}=\sum_{n\ge0}\tau^nt^n,$$
thus $R_n(\tau e_0)=\tau^n$. Since $R_n$ is a polynomial in $x_0$ and $\rho^2$, we have $\partial_jR_n=2x_j\,\dfrac{\partial R_n}{\partial(\rho^2)}=0$ at $x_1=x_2=0$, $j=1,2$.

(iii) For $t=0$, $\Phi=1$. For $t\ne0$, $1-2x_0t+r^2t^2=t^2\abs{x-t^{-1}e_0}^2$, it follows that 
$$\Phi(x,t)=\bigl(\abs t .\abs{x-t^{-1}e_0}\bigr)^{-1}$$ 
is a multiple of the Newtonian potential centered at $t^{-1}e_0$. Since $|t|<1$, $|t^{-1}e_0|=|t|^{-1}>1$, the pole lies outside $\overline{\Bu}$. Consequently, $\Phi(x,t)$ is harmonic in $\Bu$ for every $t\in(-1,1)\backslash \{0\}$. Moreover, $\Phi\in C^\infty(\Bu\times(-1,1))$ and
$$\Delta_x R_n(x) = \frac1{n!}\partial_t^n\Delta_x\Phi(x,t)\big|_{t=0}=0.$$

(iv) Let $Q=1-2x_0t+r^2t^2$. A direct computation gives $\partial_0\Phi=t(1-x_0t)Q^{-3/2}$, while
\[
\partial_t(t\Phi)=Q^{-1/2}+t(x_0-r^2t)Q^{-3/2}=Q^{-3/2}\bigl(Q+tx_0-r^2t^2\bigr)=(1-x_0t)Q^{-3/2}.
\]
Hence $\partial_0\Phi=t\,\partial_t(t\Phi)$, equivalently,
$$\sum_{n\ge0}\partial_0R_nt^n=\sum_{n\ge0}(n+1)R_nt^{n+1}.$$ 
Comparing the coefficients of $t^{n+1}$ yields $\partial_0R_{n+1}=(n+1)R_n$.
\end{proof}

We call $R_n(x)$ the zonal harmonic polynomial of degree $n$. For instance, 
\[
\aligned
& R_0=1,\qquad R_1(x)=x_0,\qquad R_2(x)=x_0^2-\frac12\rho^2,\\ 
& R_3(x) = x_0^3 - \frac32 x_0\rho^2,\qquad R_4(x)=x_0^4-3x_0^2\rho^2+\frac38\rho^4.
\endaligned
\]

\begin{lemma}[Rotation invariant; cf.\ \cite{ABR01}]\label{lem:zonalpoly}
Let $q$ be a homogeneous harmonic polynomial of degree $n$ on $\R^3$ that is invariant under all rotations about the $x_0$-axis. Then $q=c\,R_n$ for some real constant $c$.
\end{lemma}
\begin{proof}
Write 
$$q=\sum_{m=0}^nx_0^{n-m}q_m(x_1,x_2)$$ 
with $q_m$ homogeneous of degree $m$. A rotation $R_\phi$ about the $x_0$-axis acts only on $(x_1,x_2)$, hence $q_m\circ R_\phi=q_m$ for every $\phi$. Thus $q_m$ is constant on circles centred at $0$, and by homogeneity $q_m(x_1,x_2)=\gamma_m\rho^m$ for a constant $\gamma_m$. Since $q_m(-x_1,-x_2)=(-1)^mq_m(x_1,x_2)$ while $\rho^m$ is even, $\gamma_m=0$ for odd $m$. Therefore, $q=\sum_{j\ge0}\gamma_{2j}\,x_0^{\,n-2j}\rho^{2j}.$
A direct computation gives 
$$\Delta(x_0^{n-2j}\rho^{2j})=(n-2j)(n-2j-1)x_0^{n-2j-2}\rho^{2j}+(2j)^2x_0^{n-2j}\rho^{2j-2}.$$ 
Comparing the coefficient of $x_0^{n-2j}\rho^{2j-2}$ in $\Delta q=0$ for $j\ge1$ yields
\[
(n-2j+2)(n-2j+1)\gamma_{2j-2}+4j^2\gamma_{2j}=0 .
\]
Hence each $\gamma_{2j}$ is determined by $\gamma_0$, and the space of such $q$ has dimension at most $1$. By Lemma~\ref{lem:R}, it contains $R_n\ne0$. Therefore, $q = cR_n$ for some real constant $c$.
\end{proof}

\paragraph{Appell polynomials.}
For $n\ge0$, define
\begin{equation}\label{eq:An}
A_n=\frac1{n+1}\,\partial R_{n+1}=R_n-\frac1{n+1}\bigl(e_1\partial_1+e_2\partial_2\bigr)R_{n+1}.
\end{equation}
For example, $A_0=1$ and $A_1=x_0+\dfrac12(e_1x_1+e_2x_2)$.

\begin{lemma}[cf.\ \cite{BG10}]\label{lem:An}
Each $A_n$ is an $\A$-valued monogenic polynomial. Moreover, $\partial_0A_n=nA_{n-1}$ for $n\ge1$, and $A_n(te_0)=t^n$.
\end{lemma}
\begin{proof}
Since $A_n=\dfrac1{n+1}\,\partial R_{n+1}$ and $R_{n+1}$ is real-valued harmonic polynomial by Lemma~\ref{lem:R}, $A_n$ is monogenic by Lemma~\ref{lem:grad}. As the operators $\partial_0$ and $\partial$ commute, Lemma~\ref{lem:R} gives
$$\partial_0A_n=\frac1{n+1}\partial(\partial_0R_{n+1})=\partial R_n=nA_{n-1}.$$ 
On the $x_0$-axis, $\partial_1R_{n+1}=\partial_2R_{n+1}=0$ and $R_n(te_0)=t^n$, so \eqref{eq:An} yields $A_n(te_0)=t^n$.
\end{proof}

These are the zonal members of the Appell system introduced by Bock and G\"urlebeck \cite{BG10}; see also \cite{GHS08}. They play the role of the powers $z^n$ in the hypercomplex setting.

\section{Elementary bounds on the Schwarz constant}\label{sec:bounds}

Recall $\mathcal M$ and $\Cst=\sup_{\mathcal M}\abs{f'(0)}$ from the Introduction. We give two elementary bounds, based only on Section~\ref{sec:prelim}: an explicit example gives a lower bound, and an estimate from the Cauchy kernel gives an upper bound. In particular, the latter shows that $\Cst<\infty$.

\subsection{A lower bound: an explicit example}\label{sec:elem}

That the interior inequality $\abs{f(x)}\le\abs x$ fails for monogenic functions is known \cite{YQ06}; we record an explicit example with a constant, which is also a first lower bound for $\Cst$.

\begin{proposition}\label{prop:elem}
Let $f=\frac{11}{10}A_1-\frac{1}{10}A_3.$
Then $f\in\mathcal M$, $f'(0)=\frac{11}{10}$, and hence
\[
\Cst\ge\frac{11}{10}>1.
\]
Moreover,
\[
\abs{f(te_0)}=\frac{t(11-t^2)}{10}>t
\qquad\text{for }\ 0<t<1,
\]
so the inequality $\abs{f(x)}\le\abs{x}$ fails on a whole segment of the
$x_0$-axis.
\end{proposition}
\begin{proof}
By Lemma~\ref{lem:An}, $f$ is monogenic with values in $\A$.
Using 
\[
A_1 = x_0 + \frac12(e_1x_1 + e_2x_2),\qquad A_3 = x_0^3 - \frac32 x_0\rho^2 + \left(\frac32 x_0^2 - \frac38\rho^2\right)(e_1x_1 + e_2x_2),
\]
the components of $f$ are
\[
f_0=\frac{x_0(22+3\rho^2-2x_0^2)}{20},\qquad f_i=\frac{x_i(44+3\rho^2-12x_0^2)}{80}\quad (i=1,2).
\]
Hence $f(0)=0$ and $f'(0)=\partial_0f_0(0)=\frac{11}{10}$. On the $x_0$-axis, $\rho=0$ and 
$$f(te_0)=\frac{t(22-2t^2)}{20} = \frac{t(11-t^2)}{10}\, ,$$ 
so $\abs{f(te_0)}>t$ for $0<t<1$.

Let $x\in S^2$ and put $t=x_0^2\in[0,1]$. Then $\rho^2=1-t$ and a direct computation gives
\[
\abs f^2 = \frac{2209+6381t-2365t^2+175t^3}{6400}=: p(t).
\]
Since $p'(t)>0$ on $[0,1]$, $p(0) = \frac{2209}{6400}<1$, and $p(1)=1$,
we have $\abs f\leq 1$ on $S^2$. Consequently, $f\in\mathcal M$ and $\Cst\ge\abs{f'(0)}=\frac{11}{10}$.
\end{proof}

\subsection{An upper bound from the Cauchy kernel}\label{sec:cauchy}

Let $K(z)=\dfrac{\bar z}{4\pi\abs z^3}$ for $z\in\R^3\setminus\{0\}$ be the fundamental solution of $\dbar$, so that $\dbar K=\delta_0$. We use the Cauchy integral formula of Clifford analysis \cite{BDS82,GHS08,GS97}: if $g\in C^1(\overline\Bu;\Hq)$ is monogenic in $\Bu$, then
\begin{equation}\label{eq:cauchyformula}
g(x)=\int_{S^2}K(y-x)\,n(y)\,g(y)\,d\sigma(y),\qquad x\in\Bu,
\end{equation}
where $n(y)=y$ is the outward unit normal to $S^2$ at $y$.

It is well-known that for a unit vector $w\in\R^3$ and a continuous function $\varphi$ on $[-1,1]$, one has
\begin{equation}\label{eq:archimedes}
\frac1{4\pi}\int_{S^2}\varphi(\ip yw)\,d\sigma(y)=\frac12\int_{-1}^1\varphi(s)\,ds.
\end{equation}
This follows from spherical coordinates with polar axis $w$: $d\sigma=\sin\vartheta\,d\vartheta\,d\phi$ and $s=\ip yw=\cos\vartheta$.

\begin{theorem}[Cauchy-kernel bound]\label{thm:cauchy}
Let $f$ be monogenic in $\Bu$ and satisfy $\abs f\le1$. Then, for every unit vector $w\in\R^3$,
\[
\abs{\partial_wf(0)}\le C_K:=\int_0^1\sqrt{1+3s^2}\,ds=1+\frac{\ln(2+\sqrt3)}{2\sqrt3}.
\]
In particular, $\abs{f'(0)}\le C_K$ and $\norm{J_f(0)}_{\mathrm{op}}\le C_K$, so $\Cst\le C_K$.
\end{theorem}
\begin{proof}
Let $g\in C^1(\overline\Bu;\Hq)$ be monogenic in $\Bu$. For $y\in S^2$, the function $x\mapsto K(y-x)$ is smooth for $\abs x<1$, and $\partial_{w}[K(y-x)]=-(\partial_wK)(y-x)$. Differentiating \eqref{eq:cauchyformula} under the integral sign at $x=0$, we obtain
\[
\partial_wg(0)=-\int_{S^2}(\partial_wK)(y)\,n(y)\,g(y)\,d\sigma(y).
\]
Since $\partial_w\overline{z} = \overline{w}$ and $\partial_w\abs z^{-3}=-3\ip zw\,\abs z^{-5}$, we have
\[
(\partial_wK)(y)=\frac{\overline{w} - 3\ip yw\,\overline{y}}{4\pi},\qquad\text{for }\abs y=1.
\]
Using $\abs w=\abs y=1$ and $\ip{\overline{w}}{\overline{y}}=\ip wy$, we get
\[
\abs{\overline{w}-3\ip yw\overline{y}}^2=1-6\ip yw\ip{\overline{w}}{\overline{y}}+9\ip yw^2 = 1+3\ip yw^2 .
\]
Moreover, $\abs{(\partial_wK)(y)\,y\,g(y)}=\abs{(\partial_wK)(y)}\,\abs{g(y)}$. By \eqref{eq:archimedes},
\[
\begin{aligned}
\abs{\partial_wg(0)}&\le\norm g_\infty\,\frac1{4\pi}\int_{S^2}\sqrt{1+3\ip yw^2}\,d\sigma(y)\\
&=\norm g_\infty\cdot\frac12\int_{-1}^1\sqrt{1+3s^2}\,ds = \norm g_\infty\left(1+\frac{\ln(2+\sqrt3)}{2\sqrt3}\right).
\end{aligned}
\]
For $\lambda\in(0,1)$, put $g(x)=f(\lambda x)$. Then $g\in C^1(\overline\Bu)$ is monogenic in $\Bu$, $\norm g_\infty\le1$, and $\partial_wg(0)=\lambda\,\partial_wf(0)$. Thus $\lambda\abs{\partial_wf(0)}\le C_K$. Letting $\lambda\to1$, we obtain $\abs{\partial_wf(0)}\le C_K$. Taking the maximum over all unit vectors $w$ gives the corresponding bound for $\norm{J_f(0)}_{\mathrm{op}}$.
\end{proof}

\begin{remark}
In Theorem \ref{thm:cauchy}, the condition $f(0)=0$ is not needed. The upper bound is approximately $1.38017$, which is smaller than the upper bound $3/2$ for harmonic functions \cite{ABR01}. In the plane the Cauchy kernel $1/(2\pi\zeta)$ has constant modulus on the circle, which gives $\abs{f'(0)}\le 1$. In $\R^d$, the analogous computation gives 
$$\abs{\overline{w}-d\ip yw\overline{y}}^2=1+d(d-2)\ip yw^2,$$ 
which is independent of $y$ only when $d=2$; this may be the source of the excess of $C_K$ over $1$.
\end{remark}

The bounds $\frac{11}{10}\le\Cst\le C_K\approx 1.38017$ leave a gap of about $0.28$. To narrow it, we study the extremal problem directly. In the next section, using rotation covariance and a symmetrization argument, we reduce it to maximizing the linear functional $2a_2$ over zonal even functions $h$ satisfying $\abs{\partial h}\le1$ on $S^2$, a convex problem in one variable. In this form, every admissible $h$ yields a lower bound, while every admissible dual vector field yields an upper bound.

\section{The extremal problem}\label{sec:extremal}

\subsection{Symmetry and rotation covariance}\label{sec:red}

The group $SO(3)$ acts on $\R^3\cong\A$. For $R\in SO(3)$ and a real-valued function $h$ on $\Bu$, the composition $h\circ R$ is again a function on $\Bu$. By the chain rule,
\begin{equation}\label{eq:chain}
\dbar(h\circ R)(x)=R^{T}\,(\dbar h)(Rx),\qquad \Hs_{h\circ R}(x)=R^{T}\,\Hs_h(Rx)\,R .
\end{equation}

\begin{lemma}[Rotation covariance]\label{lem:rot}
Let $f=\partial h\in\mathcal M$ and $R\in SO(3)$. Then $f_R:=\partial(h\circ R)\in\mathcal M$, $\abs{f_R(x)}=\abs{f(Rx)}$ and $\Hs_{h\circ R}(0)=R^{T}\Hs_h(0)R$. Moreover, $-f\in\mathcal M$, and $\mathcal M$ is convex.
\end{lemma}
\begin{proof}
Since $\Delta$ commutes with rotations, $h\circ R$ is harmonic, and hence $f_R$ is monogenic. By \eqref{eq:chain}, 
$$\abs{f_R(x)}=\abs{\dbar(h\circ R)(x)}=\abs{R^T(\dbar h)(Rx)}=\abs{(\dbar h)(Rx)}=\abs{f(Rx)}\le1,$$
because $R^T$ is orthogonal. Moreover, $f_R(0)=\overline{R^T\overline{f(0)}}=0$. The Hessian identity is the second formula in \eqref{eq:chain} evaluated at $x=0$. 

Finally, \(\mathcal M\) is convex: if \(f_1,f_2\in\mathcal M\) and \(s\in[0,1]\), then \(sf_1+(1-s)f_2\) is again monogenic, vanishes at \(0\), and has modulus at most \(1\). Moreover, \(-f\in\mathcal M\) for every \(f\in\mathcal M\).
\end{proof}

\begin{theorem}[All measures coincide]\label{thm:meas}
For every unit vector $w\in\R^3$, the quantity 
$$c_0:=\sup_{f=\partial h\in\mathcal M}\ip{\Hs_h(0)w}{w}$$ 
is independent of $w$, and
\[
\Cst=\sup_{f\in\mathcal M}\abs{f'(0)}=\sup_{f\in\mathcal M}\norm{J_f(0)}_{\mathrm{op}}=\sup_{f\in\mathcal M}\partial_0f_0(0)=c_0<\infty .
\]
\end{theorem}
\begin{proof}
Given unit vectors $w,w'\in\R^3$, choose $R\in SO(3)$ such that $Rw'=w$. By Lemma~\ref{lem:rot}, the map $f\mapsto f_R$ is a bijection of $\mathcal M$, with inverse $f\mapsto f_{R^{-1}}$. Moreover, 
$$\ip{\Hs_{h\circ R}(0)w'}{w'}=\ip{R^T\Hs_h(0)Rw'}{w'}=\ip{\Hs_h(0)w}{w}.$$ Hence the suprema defining $c_0$ for $w$ and $w'$ agree; that is, $c_0$ is independent of $w$. In particular,
$$c_0=\sup_{f=\partial h\in\mathcal M}\ip{\Hs_h(0)e_0}{e_0}.$$

Since $-\mathcal M=\mathcal M$ and $\Hs_{-h}=-\Hs_h$, the supremum of $-\ip{\Hs_h(0)w}{w}$ over $f=\partial h\in\mathcal M$ is also $c_0$. Therefore, $\sup_{f=\partial h\in\mathcal M}\norm{\Hs_h(0)}_{\mathrm{op}}=c_0$. By Lemma~\ref{lem:grad}, 
$$\norm{J_f(0)}_{\mathrm{op}}=\norm{\Hs_h(0)}_{\mathrm{op}},\qquad \abs{f'(0)}=\abs{\Hs_h(0)e_0},\qquad \partial_0f_0(0)=\ip{\Hs_h(0)e_0}{e_0}.$$ 
Consequently,
\[
c_0=\sup_{f=\partial h\in\mathcal M}\partial_0f_0(0)\le\sup_{f=\partial h\in\mathcal M}\abs{f'(0)}=\sup_{f=\partial h\in\mathcal M}\abs{\Hs_h(0)e_0}\le\sup_{f=\partial h\in\mathcal M}\norm{\Hs_h(0)}_{\mathrm{op}}=c_0.
\]
Thus, $\sup_{f\in\mathcal M}\abs{f'(0)}=c_0$, and likewise $\sup_{f\in\mathcal M}\norm{J_f(0)}_{\mathrm{op}}=c_0$. Finally, $c_0=\Cst\le C_K$ by Theorem~\ref{thm:cauchy}.
\end{proof}

\subsection{Zonal reduction}\label{sec:zonal}

\begin{definition}
A function $h$ on $\Bu$ is zonal if $h\circ R=h$ for every rotation $R$ about the $x_0$-axis, and even if $h(-x)=h(x)$. Let $\mathcal Z\subset\mathcal M$ denote the set of functions $f=\partial h$ with $h$ zonal and even.
\end{definition}

\begin{lemma}[Expansion of zonal even functions]\label{lem:expansion}
Let $h$ be a zonal even harmonic function on $\Bu$. Then there exist real numbers $a_0,a_2,a_4,\dots$ such that
\[
h=a_0+\sum_{k\ge1}a_{2k}\,R_{2k},
\]
where the series converges locally uniformly in $\Bu$, together with all its derivatives. Consequently,
$$\Hs_h(0)=a_2\diag(2,-1,-1),\qquad J_f(0) = 2a_2\diag\left(1,\frac12,\frac12\right)\quad \text{for } f=\partial h,$$
and 
$$f=\partial h=\sum_{k\ge1}2k\,a_{2k}A_{2k-1}.$$
\end{lemma}
\begin{proof}
A harmonic function on the ball $\Bu$ admits a locally uniformly convergent expansion $h=\sum_{n\ge0}h_n$, where each $h_n$ is a homogeneous harmonic polynomial of degree $n$ (see \cite{ABR01}); the expansion may be differentiated term by term. 

Since $h$ is invariant under rotations about the $x_0$-axis, and rotations preserve homogeneity, each $h_n$ is invariant as well, by uniqueness of the expansion. Similarly, the degree-$n$ component of $h(-\,\cdot)$ is $(-1)^nh_n$, thus $h_n=0$ for odd $n$. By Lemma~\ref{lem:zonalpoly}, $h_n=a_nR_n$ for even $n\ge2$, while $h_0=a_0$ is constant. 

Since $R_2=x_0^2-\frac12\rho^2$, we have $\Hs_{R_2}=\diag(2,-1,-1)$. Moreover, $\Hs_{R_{2k}}(0)=0$ for $k\ge2$, so $\Hs_h(0)=a_2\diag(2,-1,-1)$. By Lemma~\ref{lem:grad}, $J_f(0)=C\,\Hs_h(0)=a_2\diag(2,1,1)$, where $C=\diag(1,-1,-1)$. Finally, by \eqref{eq:An}, 
$$f=\partial h=\sum_{k\ge1} a_{2k}\partial R_{2k}=\sum_{k\ge1} 2k\,a_{2k}A_{2k-1}.$$
\end{proof}

\begin{lemma}[Attainment]\label{lem:attain}
There exists $f^*\in\mathcal M$ such that $\partial_0f^*_0(0)=\abs{(f^{*})'(0)}=\Cst$. Such a function $f^*$ is called an extremal.
\end{lemma}
\begin{proof}
By Theorem~\ref{thm:meas}, $\Cst=\sup_{f\in\mathcal M}\partial_0f_0(0)<\infty$; so there exist $f_k\in\mathcal M$ with $\partial_0(f_k)_0(0)\to\Cst$. Since $\abs{f_k}\le1$, the components of $f_k$ are harmonic and bounded by $1$; hence, by the normal family theorem for harmonic functions (see \cite{ABR01}), a subsequence converges locally uniformly in $\Bu$, together with all derivatives, to some $f^*$. Since a locally uniform limit of monogenic functions is monogenic (see \cite{GHS08}), $f^*$ is monogenic, and the conditions $f^*(0)=0$ and $\abs{f^*}\le1$ are preserved; thus $f^*\in\mathcal M$. Since $\partial_0f_0(0)$ involves only first derivatives at $0$, we get 
$$\partial_0f^*_0(0)=\lim_{k\to\infty}\partial_0(f_k)_0(0)=\Cst.$$
Consequently, 
$$\Cst=\partial_0f^*_0(0)\le\abs{(f^*)'(0)}\le\Cst,$$ 
which gives the claim.
\end{proof}

\begin{theorem}[Zonal reduction]\label{thm:zonal}
For every $f=\partial h\in\mathcal M$, its average $g:=\partial\widehat{h}$, where
$$\widehat{h}(x)=\frac1{4\pi}\int_0^{2\pi}\big(h(R_\phi x)+h(-R_\phi x)\big)\,d\phi$$ 
and $R_\phi$ is the rotation about the $x_0$-axis by the angle $\phi$, belongs to $\mathcal Z$ and satisfies $\partial_0g_0(0)=\partial_0f_0(0)$. Consequently,
\[
\Cst=\sup\{2a_2(h):\ f=\partial h\in\mathcal Z\},
\]
and extremal sequences may be chosen to be zonal. Moreover, if $f^*\in\mathcal M$ is extremal (Lemma~\ref{lem:attain}), then its average $g\in\mathcal Z$ satisfies $\partial_0g_0(0)=\Cst$; hence $g$ is again extremal, and $J_g(0)=\Cst\,\diag(1,\frac12,\frac12)$.
\end{theorem}
\begin{proof}
Let $f=\partial h\in\mathcal M$ and let $R_\phi$ denote the rotation by $\phi$ about the $x_0$-axis. Put 
$$\widetilde h=\frac1{2\pi}\int_0^{2\pi}h\circ R_\phi\,d\phi.$$ 
The integral converges locally uniformly, together with all derivatives, since harmonic functions are smooth; hence $\widetilde h$ is harmonic and 
$$\partial\widetilde h=\frac1{2\pi}\int_0^{2\pi}f_{R_\phi}\,d\phi.$$
Each $f_{R_\phi}$ lies in $\mathcal M$ by Lemma~\ref{lem:rot}. Thus $\partial\widetilde h$ is monogenic, vanishes at $0$ and has modulus at most $1$, so $\partial\widetilde h\in\mathcal M$. The function $\widetilde h$ is zonal, and, since $R_\phi e_0=e_0$, \eqref{eq:chain} gives 
$$\ip{\Hs_{\widetilde h}(0)e_0}{e_0}=\frac1{2\pi}\int_0^{2\pi}\ip{R_\phi^T\Hs_h(0)R_\phi e_0}{e_0}d\phi=\ip{\Hs_h(0)e_0}{e_0}.$$

Next put $h^-(x)=h(-x)$. Then $\partial h^-(x)=-f(-x)$, so $\partial h^-\in\mathcal M$, and $\Hs_{h^-}(0)=\Hs_h(0)$. The function $\widehat h:=\frac12(\widetilde h+\widetilde h^-)$ is zonal and even, $\partial\widehat h\in\mathcal M$ by convexity and 
$$\ip{\Hs_{\widehat h}(0)e_0}{e_0}=\ip{\Hs_h(0)e_0}{e_0}.$$ 
This gives $g=\partial\widehat h\in\mathcal Z$. 
By Lemma~\ref{lem:expansion}, $g=\sum_{k\ge1}2k\,a_{2k}A_{2k-1}$. Thus, we have
$$\partial_0f_0(0) = \ip{\Hs_h(0)e_0}{e_0} = \ip{\Hs_{\widehat h}(0)e_0}{e_0} = \partial_0g_0(0) = 2a_2.$$ 
By Theorem~\ref{thm:meas}, 
$$\Cst=\sup_{f\in\mathcal M}\partial_0f_0(0)=\sup\{2a_2:\ g=\partial  \widehat h\in\mathcal Z\}.$$ 
Finally, if $f^*$ is extremal, i.e.\ $\partial_0f^*_0(0)=\Cst$, then its average $g$ satisfies $2a_2=\Cst$, hence $g$ is extremal and $J_g(0)=2a_2\diag(1,\frac12,\frac12)=\Cst\diag(1,\frac12,\frac12)$.
\end{proof}

\subsection{Sharper bounds from the zonal problem}\label{sec:duality}

\subsubsection*{Lower bounds from polynomials}

Let $\mathcal H_N$ denote the span of $R_2,R_4,\dots,R_N$, where $N$ is even.

\begin{proposition}[Primal restriction]\label{prop:primal}
If $h\in\mathcal H_N$ is non-zero and $\Gamma\ge\max_{S^2}\abs{\partial h}^2$, then 
\[
\Cst\ge \frac{2a_2(h)}{\sqrt\Gamma}.
\] 
The optimal values 
\[
V_N = \sup\left\{2a_2(h):\ h\in\mathcal H_N,\ \abs{\dbar h}\le1\text{ on }S^2\right\}
\] 
are non-decreasing in the even integer $N$ and converge to $\Cst$.
\end{proposition}
\begin{proof}
The function $f=\partial h/\sqrt\Gamma$ is monogenic and satisfies $\abs f\le1$ on $S^2$, hence $\abs f\le1$ in $\Bu$ by Lemma~\ref{lem:max}. Moreover, $f(0)=0$ because $h$ has no terms of degree less than $2$. Thus $f\in\mathcal M$ and 
\[
\Cst\ge \partial_0f_0(0) = \frac{2a_2(h)}{\sqrt\Gamma}.
\] 
The monotonicity of $V_N$ follows from the inclusion $\mathcal H_N\subset\mathcal H_{N+2}$. Note that $V_N\le\Cst$, since every admissible $h$ gives $\partial h\in\mathcal Z$.

For convergence, let $g=\partial u\in\mathcal Z$, where $u=\sum_{n\ge2,n\text{ even}} a_{n}R_{n}$ as in Lemma~\ref{lem:expansion}, and let $0<\lambda<1$. Define
\[
u_\lambda(x) = \frac{u(\lambda x)}{\lambda}.
\] 
Then $\partial u_\lambda(x)=g(\lambda x)$, so $\abs{\dbar u_\lambda}=\abs{\partial u_\lambda}\le1$ on $\overline\Bu$. Its coefficients are $\lambda^{n-1}a_n$, and it is harmonic on a ball of radius $1/\lambda>1$. Let $u_\lambda^N\in\mathcal H_N$ for even $N\ge2$ denote the Taylor partial sum of $u_\lambda$ of degree $N$. Because $u_\lambda^N$ and $u_\lambda$ have the same homogeneous component of degree $2$, we have $a_2(u^N_\lambda)=\lambda a_2$. Moreover, $u_\lambda^N\to u_\lambda$ in $C^1(\overline\Bu)$, and $\abs{\dbar u_\lambda}\le1$, thus
\[
\limsup_{N\to\infty}\max_{S^2}\abs{\dbar u_\lambda^N}\le1.
\] 
Put $\epsilon_N:=\max_{S^2}\abs{\dbar u_\lambda^N}-1$ and $\epsilon_N^+:=\max(\epsilon_N,0)$, so $\epsilon_N^+\to0$. Defining $\widetilde u_\lambda^N:=u_\lambda^N/(1+\epsilon_N^+)$ gives $\max_{S^2}\abs{\dbar\widetilde u_\lambda^N}\le1$. Also $a_2(\widetilde u_\lambda^N)=a_2(u^N_\lambda)/(1+\epsilon_N^+)=\lambda a_2/(1+\epsilon_N^+)$; hence
\[
V_N\ge \frac{2\lambda a_2}{1+\epsilon_N^+}.
\]
Since $\epsilon_N^+\to0$, for every $0<\lambda<1$,
$$\liminf_{N\to\infty} V_N\ge2\lambda a_2.$$ 
Letting $\lambda\uparrow1$ gives $\liminf_{N\to\infty} V_N\ge2a_2$. Taking the supremum over $g\in\mathcal Z$ and using Theorem~\ref{thm:zonal}, we get 
\[
\liminf_{N\to\infty} V_N\ge\Cst.
\]
Together with $V_N\le\Cst$ and the monotonicity, this gives $V_N\to\Cst$.
\end{proof}

Proposition~\ref{prop:primal} turns every $h\in\mathcal H_N$ into a lower bound for $\Cst$. The quality of the bound depends on the choice of $h$: one wants $a_2(h)$ large while $\abs{\dbar h}$ stays small on $S^2$.

Let $(r,\vartheta,\varphi)$ be spherical coordinates with polar axis $x_0$, that is, 
\[
x_0=r\cos\vartheta,\qquad x_1=r\sin\vartheta\cos\varphi,\qquad x_2=r\sin\vartheta\sin\varphi.
\] 
Put $\omega=\cos\varphi\,e_1+\sin\varphi\,e_2$, and
\[
\mathbf e_r=\cos\vartheta+\omega\sin\vartheta,\qquad
\mathbf e_\vartheta=-\sin\vartheta+\omega\cos\vartheta,\qquad
\mathbf e_\varphi=-\sin\varphi\,e_1+\cos\varphi\,e_2 .
\]
Under $\mathcal A\cong\R^3$, these are the orthonormal vectors of the spherical frame, and the Cauchy-Riemann operator has the form
\[
\dbar=\mathbf e_r\,\partial_r+\mathbf e_\vartheta\,\frac1r\,\partial_\vartheta+\mathbf e_\varphi\,\frac1{r\sin\vartheta}\,\partial_\varphi .
\]

For $h\in\mathcal H_N$, since $h=\sum_{n\ge2,n\text{ even}}a_nr^nP_n(\cos\vartheta)$ depends only on $r$ and $\vartheta$,
\[
\dbar h = \mathbf e_r\,\partial_rh + \mathbf e_\vartheta\,\frac1r\,\partial_\vartheta h .
\]
Moreover, $\partial_\vartheta[P_n(\cos\vartheta)]=-\sin\vartheta\,P_n'(\cos\vartheta)$. Hence on $S^2$, with $c=\cos\vartheta$,
\begin{equation}\label{eq:G}
\abs{\dbar h}^2=\abs{\partial h}^2=G(c):=\Bigl(\sum_n n\,a_nP_n(c)\Bigr)^2+(1-c^2)\Bigl(\sum_na_nP_n'(c)\Bigr)^2 .
\end{equation}
For a given $h$, the smallest admissible constant in Proposition~\ref{prop:primal} is $\Gamma=\max_{c\in[-1,1]}G(c)$, and we obtain the lower bound
\[
\Cst\ge\frac{2a_2(h)}{\sqrt{\max_{c}G(c)}} .
\]
The best lower bound for $\Cst$ obtainable from $\mathcal H_N$ is exactly
\[
V_N=\max\{2a_2(h):\ h\in\mathcal H_N,\ G(c)\le1\text{ for all }c\in[-1,1]\}.
\]
The space $\mathcal H_N$ has dimension $N/2$, so this is a search over finitely many parameters $a_2,a_4,\dots,a_N$. Moreover, the problem is convex, since the objective $2a_2$ is linear and the constraint $\abs{\partial h}\le1$ is convex in the coefficients. There exist efficient algorithms that find its global optimum $V_N$. We do not pursue the numerical solution here, but emphasize that the best lower bound for $\Cst$ obtainable in $\mathcal H_N$ can indeed be found. By Proposition~\ref{prop:primal}, $V_N$ increases to $\Cst$ as $N\to\infty$.

\begin{example}\label{ex:lower}
We illustrate this with the cases $N=2,4,6$. 

Consider $N=2$ and $h=\frac12R_2$. By~\eqref{eq:G}, with $t=c^2$, 
\[
G=\frac{1+3t}{4}.
\] 
Then $\max_{t\in[0,1]}G(t) = G(1) = 1$. Put $\Gamma = 1$; hence $\Cst\ \ge\ 2\cdot\frac12/\sqrt{1}=1$, which is the trivial bound.
For $N=4$ and $N=6$, the same computation gives Table~\ref{tab:lower}.

\begin{table}[ht]
\centering
\renewcommand{\arraystretch}{1.5}
\begin{tabular}{|p{0.15\linewidth}|p{0.30\linewidth}|p{0.40\linewidth}|}
\hline
 & $N=4$ & $N=6$ \\
\hline
$h_N$ & $\frac58 R_2-\frac3{32}R_4$
      & $\frac58 R_2-\frac7{64}R_4+\frac3{128}R_6$ \\
\hline
$a_2$ & $\frac58$ & $\frac58$ \\
\hline
$G(t)$ & $(2401+6645t-7485t^2+1575t^3)/4096$
       & $(727609+1701855t-4295830t^2+4152526t^3-1726515t^4+392931t^5)/1048576$ \\
\hline
$\max G$ & $0.98995\ldots$ (at $t_0\approx 0.5338$)
         & $(61/64)^2=0.90844\ldots$ (at $t=1$) \\
\hline
$\Gamma$ & $0.99$ & $(61/64)^2$ \\
\hline
$2a_2/\sqrt{\Gamma}$ & $1.25629\ldots$ & $80/61=1.31147\ldots$ \\
\hline
\end{tabular}
\caption{Lower bounds for $\Cst$ from $h_N\in\mathcal H_N$; here $t=c^2$ and $\Gamma\ge\max_{[0,1]}G$.}
\label{tab:lower}
\end{table}
\end{example}

\subsubsection*{Upper bounds from dual fields}

The lower bounds use a single zonal harmonic function. For upper bounds, we use two zonal functions, which generate the radial and tangential components of a vector field on the sphere.

For $x\ne0$ put $y=x/\abs x\in S^2$, so that $\mathbf e_r=y$. For a real function $g$ of class $C^1$ near $S^2$, setting $r=1$ in the formula for $\dbar$ gives
\[
\dbar g(y)=(\partial_rg)(y)\,y+\dbar_Sg(y),\qquad
\dbar_Sg:=\mathbf e_\vartheta\,\partial_\vartheta g+\mathbf e_\varphi\,\frac1{\sin\vartheta}\,\partial_\varphi g ,
\]
where $\partial_rg=\ip{\dbar g}{y}$ is the radial part and $\dbar_Sg$ is the tangential part, $\ip{\dbar_Sg}{y}=0$. Since $\dbar_Sg$ involves only derivatives along $S^2$, it depends only on $g|_{S^2}$. Note that, for a function $\kappa(c)$ of $c=y_0=\cos\vartheta$, we have $\partial_\varphi\kappa=0$ and $\partial_\vartheta\kappa=-\sin\vartheta\,\kappa'(c)$; hence $\abs{\dbar_S\kappa}^2=(1-c^2)\kappa'(c)^2$. 

Let $u=\sum_{n\ge0}\alpha_nP_n(c)$ and $v=\sum_{n\ge1}\beta_nP_n(c)$ with $\alpha_n,\beta_n=O(q^n)$ for some $0<q<1$; then $u,v$ are real analytic on $[-1,1]$ (see \cite{Sz75}). Let $W=u\,y+\dbar_Sv$ on $S^2$, that is, $W$ has radial component $u$ and tangential component $\dbar_Sv$. Since these components are orthogonal,
\[
\abs W^2=u^2+(1-c^2)\,v'(c)^2 .
\]

The Green formula for $g_1,g_2\in C^2(S^2)$ reads 
\[
\int_{S^2}\ip{\dbar_Sg_1}{\dbar_Sg_2}\,d\sigma=-\int_{S^2}g_1\,\Delta_Sg_2\,d\sigma ,
\] 
where $\Delta_S$ is the Laplace--Beltrami operator on $S^2$. Moreover, by Legendre's equation, $\Delta_S[P_n(c)]=-n(n+1)\,P_n(c)$.

\begin{lemma}[Duality identity]\label{lem:dual}
Let $h=\sum_{n\ge2}a_nR_n$ be harmonic in a neighbourhood of $\overline\Bu$, and let $W=u\,y+\dbar_Sv$ be as above, with $u=\sum_{n\ge0}\alpha_nP_n(c)$, $v=\sum_{n\ge1}\beta_nP_n(c)$, $\alpha_n,\beta_n=O(q^n)$, $0<q<1$. Then
\[
\int_{S^2}\ip{\dbar h}{W}\,d\sigma=\sum_{n\ge2}\bigl(n\alpha_n+n(n+1)\beta_n\bigr)\frac{4\pi}{2n+1}\,a_n .
\]
\end{lemma}
\begin{proof}
On $S^2$, $\dbar h=\partial_rh\,y+\dbar_Sh$ and $W=u\,y+\dbar_Sv$ with $\dbar_S(\cdot)\perp y$, so 
\[
\ip{\dbar h}W=u\,\partial_rh+\ip{\dbar_Sh}{\dbar_Sv}.
\]
Moreover, on $S^2$, $h=\sum_{n\ge2}a_nP_n(c)$ and $\partial_rh=\sum_{n\ge2}na_nP_n(c)$, and both series converge absolutely and uniformly, by the interior estimates for harmonic functions (see \cite{ABR01}), since $h$ is harmonic in a ball of radius larger than $1$. By the Green formula, 
\[
\int_{S^2}\ip{\dbar_Sh}{\dbar_Sv}\,d\sigma=-\int_{S^2} h\,\Delta_Sv\,d\sigma,
\] 
where $-\Delta_Sv=\sum_{n\ge1} n(n+1)\beta_nP_n$ converges uniformly because $\beta_n=O(q^n)$. Now 
\[
\int_{S^2}P_nP_m\,d\sigma=2\pi\int_{-1}^1P_nP_m\,dc=\frac{4\pi}{2n+1}\delta_{nm} ,
\] 
and integrating term by term,
\[
\int_{S^2} u\,\partial_rh\,d\sigma=\sum_{n\ge2} n\alpha_na_n\frac{4\pi}{2n+1} ,\qquad -\int_{S^2} h\,\Delta_Sv\,d\sigma=\sum_{n\ge2} n(n+1)\beta_na_n\frac{4\pi}{2n+1}.
\]
Adding the two integrals proves the claim.
\end{proof}

\begin{proposition}[Dual restriction]\label{prop:dual}
Let $\alpha_0\in\mathbb R$ and $\beta_n\in\mathbb R$ ($n\ge2$ even), with finitely many $\beta_n\ne0$, and put
$\alpha_2=\frac5{4\pi}-3\beta_2$, $\alpha_n=-(n+1)\beta_n$ ($n\ge4$ even), and $\alpha_n=\beta_n=0$ for odd $n$.
Let $u,\,v,\,W$ be as above. Then
\[
\Cst\le\int_{S^2}\abs W\,d\sigma=4\pi\int_0^1\sqrt{u^2+(1-c^2)\,v'(c)^2}\,dc .
\]
\end{proposition}
\begin{proof}
Let $f=\partial h\in\mathcal Z$ and $0<\lambda<1$. Then $h_\lambda(x) = h(\lambda x)/\lambda$ is harmonic near $\overline\Bu$, zonal, $\dbar h_\lambda(0)=\overline{f(0)}=0$, $\abs{\dbar h_\lambda(y)}=\abs{f(\lambda y)}\le1$ for $y\in S^2$, and its coefficients are $\lambda^{n-1}a_n$ (Proof of Proposition~\ref{prop:primal}). 

For $n\ge4$, $n\alpha_n+n(n+1)\beta_n=0$ by the choice of $\alpha_n$, and $2\alpha_2+6\beta_2=\frac{10}{4\pi}$, so the coefficient of $a_2$ in Lemma~\ref{lem:dual} is $\frac{10}{4\pi}\cdot\frac{4\pi}5=2$. Therefore
\[
2\lambda a_2=\int_{S^2}\ip{\dbar h_\lambda}{W}\,d\sigma\le\int_{S^2}\abs{\dbar h_\lambda}\,\abs W\,d\sigma \le \int_{S^2}\abs W\,d\sigma.
\]
Letting $\lambda\to1$ gives $2a_2\le\int_{S^2}\abs W\,d\sigma$ for every $f\in\mathcal Z$, and the supremum over $f$ is $\Cst$ by Theorem~\ref{thm:zonal}. The last form of the integral holds because $u,v$ are even in $c$ and $d\sigma=2\pi\,dc$.
\end{proof}

\begin{remark}
Lemma~\ref{lem:dual} holds for arbitrary $\alpha_n,\beta_n$ of the above type; the choice made here only serves to make the right-hand side equal to $2a_2$. Other choices give upper bounds for other linear functionals of $(a_n)$, by the same argument.
\end{remark}

\begin{corollary}[A sharper upper bound]\label{cor:U}
In Proposition~\ref{prop:dual}, let 
\[
\alpha_0=\frac{a}{4\pi},\qquad \beta_2=\frac{1}{4\pi},\qquad \alpha_2=\frac5{4\pi}-3\beta_2=\frac2{4\pi},\qquad \alpha_n=\beta_n=0\ (n\ge3)
\] 
for $0<a<1$. Then
\[
\Cst\ \le\ U(a):=\frac12\Bigl[a+2+\frac{(1-a)^2}{\sqrt{6a+3}}\,\operatorname{arsinh}\frac{\sqrt{6a+3}}{1-a}\Bigr].
\]
\end{corollary}
\begin{proof}
Since $P_2=\frac12(3c^2-1)$ and $P_2'=3c$, 
\[
4\pi u=a+(3c^2-1)=3c^2-1+a,\qquad 4\pi v'(c)=3c.
\] 
Hence
\[
(4\pi)^2\bigl(u^2+(1-c^2)v'(c)^2\bigr)=(3c^2-1+a)^2+(1-c^2)(3c)^2=(1-a)^2+(6a+3)c^2.
\]
Thus, by Proposition~\ref{prop:dual},
\begin{align*}
\Cst & \le4\pi\int_0^1\sqrt{u^2+(1-c^2)\,v'(c)^2}\,dc \\
& =\int_0^1\sqrt{(1-a)^2+(6a+3)c^2}\,dc=\frac12\Bigl[a+2+\frac{(1-a)^2}{\sqrt{6a+3}}\,\operatorname{arsinh}\frac{\sqrt{6a+3}}{1-a}\Bigr].
\end{align*} 
\end{proof}

\begin{example}\label{ex:upper}
For $a=0.31$ we have 
\[
a+2=2.31,\qquad \frac{(1-a)^2}{\sqrt{6a+3}}=0.21596\ldots,\qquad \operatorname{arsinh}\frac{\sqrt{6a+3}}{1-a}=1.87836\ldots
\]
Rounding up all terms gives
\[
\Cst\le U(0.31) < \frac12(2.31+0.21597\cdot1.87837)\ <\ 1.35784 .
\]
\end{example}

\begin{remark}
Together with Example~\ref{ex:lower},
\[
\frac{80}{61}\ \le\ \Cst\ \le\ 1.35784 ,
\]
so $\Cst$ is located within $0.0464$ by hand-checkable estimates. The gap comes mostly from the lower bound: a polynomial with three parameters cannot make $\abs{\dbar h}$ close to $1$ on all of $S^2$. Higher-degree polynomials on the primal side and richer dual fields on the dual side narrow the interval, but they only ever give two-sided bounds, not a closed form for $\Cst$.
\end{remark}

\begin{remark}
The components of $f\in\mathcal M$ are harmonic with $\abs f\le1$, so one may ask the same question for all harmonic maps $\Bu\to\overline\Bu$ fixing $0$. For a real harmonic $u$ on $\Bu$ with $\abs u\le1$, one has $\abs{\nabla u(0)}\le\frac32$, and the constant is sharp \cite{ABR01}. If $g$ is a harmonic map $\Bu\to\overline\Bu$ fixing $0$, and $\xi$ is an arbitrary unit vector, this estimate applied to $u=\ip g\xi$ gives 
\[
\ip{\partial_0g(0)}{\xi}=\partial_0u(0)\le\abs{\nabla u(0)}\le\frac32;
\] 
taking the supremum over $\xi$ yields $\abs{\partial_0g(0)}\le\frac32$. This is sharp as well: let $u$ be the Poisson integral of $\operatorname{sign}y_0$; then $u(0)=0$, $\abs u\le1$, and $\partial_0u(0)=\frac3{4\pi}\int_{S^2}\abs{y_0}\,d\sigma=\frac32$; hence $g=u\, e_0$ satisfies $\abs{\partial_0g(0)}=\frac32$. Monogenicity therefore lowers the sharp constant: $\Cst\le C_K<\frac32$ already, and Section~\ref{sec:duality} gives $\Cst\le 1.35784$.
\end{remark}

\section{Conclusion}\label{sec:concl}

We determine bounds on the Schwarz constant $\Cst$ of the hypercomplex derivative for monogenic functions $f\colon\Bu\subset\R^3\to\A$ satisfying the classical conditions $f(0)=0$ and $\abs f\le1$.
First, an explicit example shows that the classical Schwarz lemma fails in this setting, that is, $\Cst>1$; this is also observed in \cite{YQ06}. Elementary estimates then give $\frac{11}{10}\le\Cst\le C_K=1+\frac{\ln(2+\sqrt3)}{2\sqrt3}$.

To sharpen these bounds we use the structure of the extremal problem. By rotation covariance, the natural measures of the size of the derivative at $0$ give the same constant, and the problem reduces to the class of zonal functions, where it becomes a convex one-dimensional problem with a dual. Restricting the primal problem to polynomials, that is, optimizing over one zonal harmonic function $h$ with $\abs{\dbar h}\le1$ on $S^2$, gives lower bounds; restricting the dual fields to those of the form $W=u\,y+\dbar_Sv$, that is, optimizing over two zonal functions $u,v$, gives upper bounds. In this way we obtain
\[
\frac{80}{61}\ \le\ \Cst\ \le\ 1.35784 .
\]

We outline several directions for future work. The most immediate concerns $\Cst$ itself: a closed form, or at least a description of the extremal function, including its existence and boundary behaviour. A second question is whether $\abs{f(x)}\le\Cst\abs x$ holds on all of $\Bu$ for $f\in\mathcal M$, which would give a Schwarz lemma with constant $\Cst$. Beyond this, one may look for the sharp form of a Schwarz--Pick estimate, bounding $\abs{f'(x)}$ in terms of $1-\abs{f(x)}^2$ and $1-\abs x^2$. More broadly, we would like to explore how far the methods of Schur analysis extend to monogenic functions, with the Schwarz inequality studied here as a first step.

\section*{Statements and Declarations}
\paragraph{Competing interests.} The authors have no competing interests to declare that are relevant to the content of this article.
\paragraph{ORCID.} Hung Manh Nguyen: \url{https://orcid.org/0009-0004-6503-6693}

\end{document}